\documentclass[a4paper,12pt]{article}
\title{
\LARGE{
On the natural extension of a map 
with a Siegel or Cremer point
}
}
\author{
Carlos Cabrera and 
Tomoki Kawahira
\thanks{
This research was partially supported by FY2010 Researcher
Exchange Program between JSPS and CONACYT.
}
}

\usepackage{enumerate}
\usepackage{graphicx}
\usepackage{amsmath}
\usepackage{amssymb}
\usepackage{amsfonts}
\usepackage{theorem}

\newtheorem{thm}{Theorem}[section]
\newtheorem{prop}[thm]{Proposition}
\newtheorem{lem}[thm]{Lemma}

\theorembodyfont{\rmfamily}
\newtheorem{pf}{Proof.}

\newcommand{\norm}[1]{{\left\| #1 \right\|}}

\newcommand{\paren}[1]{{\left( #1 \right)}}
\newcommand{\braces}[1]{{\left\{ #1 \right\}}}
\newcommand{\C}{\mathbb{C}}
\newcommand{\Cbar}{\overline{\mathbb{C}}}

\newcommand{\D}{\mathbb{D}}
\newcommand{\N}{\mathbb{N}}
\newcommand{\al}{{\alpha}}
\newcommand{\e}{\epsilon}
\newcommand{\cR}{{\mathcal{R}}}
\newcommand{\cN}{{\mathcal{N}}}

\newcommand{\zhat}{\hat{z}}
\newcommand{\fhat}{\hat{f}}
\newcommand{\Hhat}{\widehat{H}}
\newcommand{\Deltahat}{\widehat{\Delta}}
\newcommand{\id}{\mathrm{id}}
\newcommand{\dens}{\mathrm{dens}}
\newcommand{\st}{\,:\,}
\newcommand{\QED}{\hfill $\blacksquare$}

\newcommand{\parag}[1]{\par\medskip\noindent{\bf #1}}

\begin{document}
\maketitle

\begin{abstract}
In this note we show that the regular part of 
the natural extension (in the sense of Lyubich and Minsky \cite{LM}) 
of quadratic map $f(z) = e^{2 \pi i \theta}z + z^2$ with irrational $\theta$ of bounded type has only parabolic leaves except the invariant lift of the Siegel disk. 


We also show that though the natural extension of a rational function with a Cremer fixed point has a continuum of irregular points, 
it can not supply enough singularity to apply 
the Gross star theorem to find hyperbolic leaves.
\end{abstract}

\section{Introduction}
\parag{Natural extension and regular part.}
Let $f:\Cbar \to \Cbar = \C \cup \{\infty\}$ be a rational function of degree $\ge 2$.
It generates a non-invertible dynamical system $(f,\Cbar)$
but it also generates an invertible dynamics
in the space of ``backward orbits" (the inverse limit) 
$$
\cN_f :=  
\braces{\zhat = (z_{-n})_{n \ge 0} 
\st z_0 \in \Cbar,~z_{-n} = f(z_{-n-1})}
$$
with action
$$
\fhat((z_0, z_{-1}, \ldots )) :=
(f(z_0), f(z_{-1}), \ldots ) = 
(f(z_0), z_0, z_{-1}, \ldots ).
$$
We say $\cN_f$ (with dynamics by $\fhat$) is the \textit{natural extension} of $f$, with topology induced by $\Cbar\times\Cbar\times \cdots$.
Note that the {\it projection} $\pi:\cN_f \to \Cbar$ defined by $\pi(\zhat):= z_0$ 
semiconjugates $\fhat$ and $f$.

In 1990's, M.Lyubich and Y.Minsky \cite{LM} introduced 
the theory of hyperbolic 3-lamination associated with rational functions, which is analogous to the theory of hyperbolic 3-manifold 
for Kleinian groups. 
The theory is based on the study of the natural extension, 
in particular the subset called the \textit{regular part}
(or \textit{regular leaf space}), defined as follows:
The point $\zhat =(z_0, z_{-1}, \ldots)$ is \textit{regular}
if there exists a neighborhood $U_0$ of $z_0$ whose pull-back
$\cdots \to U_{-1} \to U_0$ along $\zhat$ (i.e., $U_{-n}$ is the connected component of $f^{-1}(U_{n-1})$ containing $z_{-n}$) is eventually univalent. 
The \textit{regular part} $\cR_f$ of $\cN_f$ is the set of all regular points, and we say the point in $\cN_f-\cR_f$ is \textit{irregular}.
The regular part is invariant under $\fhat$,
and each path-connected component (``leaf") of the regular part
possesses a Riemann surface structure isomorphic to $\C$,
$\D$, or an annulus. 
(The annulus appears only when $f$ has a Herman ring.)

\parag{Type problem.}
When the critical orbits of $f$ behave nicely,
we may regard $\cR_f$ as a Riemann surface lamination
with all leaves isomorphic to $\C$.
Such a situation yields some nice properties of dynamics, 
like rigidity, or existence of conformal invariant measures
on the lamination.  
For example, this is the case when $f$ has no recurrent critical points in the Julia set \cite[Prop.4.5]{LM}. 
Another intriguing case is when $f$ is an infinitely renormalizable quadratic map with a persistently recurrent critical point \cite[Lem.3.18]{KL}.

For general cases, it is questioned in \cite[\S 4, \S 10]{LM} when we have leaves of \textit{hyperbolic type}, 
especially leaves isomorphic to $\D$. 
(The counterpart, leaves isomorphic to $\C$, are conventionally called \textit{parabolic}.)
Easy examples of hyperbolic leaves are provided by
invariant lifts of Siegel disks and Herman rings.
Non-rotational hyperbolic leaves (that are rather non-trivial) 
are found by J.Kahn and by J.Rivera-Letelier. 
Readers may find details in the paper by J.Kahn, M.Lyubich, and L.Rempe \cite[\S 3]{KLR}, that can be summarized 
as follows: When the Julia set is contained in the postcritical set,
there are infinitely many hyperbolic leaves. 
In this case such a leaf may not touch the (lifted) Julia set
 in the natural extension.
However, by using 
the tuning technique and
\textit{the Gross star theorem},
which gives a necessity condition of a leaf to be parabolic,
we may construct a quadratic polynomial with hyperbolic leaves
that do intersect with the Julia set. 
In both cases, a recurrent critical point plays a crucial role.

In the quest of new non-rotational hyperbolic leaves,
it is natural to ask the following question:
{\it Is there any non-rotational hyperbolic leaf when $f$ has 
an irrationally indifferent fixed point?} 
Because existence of such a fixed point implies 
existence of a recurrent critical point
whose postcritical set is a continuum, 
and it seems really close to the situations in \cite{KLR}.
The aim of this note is to give some results 
on this question. 

\parag{Siegel disk of bounded type.}
$f(z) = e^{2 \pi i \theta}z + z^2$ with irrational $\theta$ of 
 bounded type has a Siegel disk $\Delta$ centered at the origin,
 whose boundary $\partial \Delta$ is a quasicircle.

In Section 2 of this note we give a proof of 

\begin{thm}[No hyperbolic leaf except the Siegel disk]\label{thm_siegel}
In the regular part of the natural extension $f:\Cbar \to \Cbar$,
the only hyperbolic leaf is the invariant lift $\Deltahat$ of the Siegel disk.
\end{thm}

In the proof we use Lyubich and Minsky's criteria
for parabolic leaves,
\textit{uniform deepness} of the postcritical set, 
and one of McMullen's results on bounded type Siegel disks.

\parag{Feigenbaum maps.}
It would be worth mentioning that 
the same method as the proof of Theorem \ref{thm_siegel} 
can be applied to a class of 
infinitely renormalizable quadratic maps,
called \textit{Feigenbaum maps}
(Precise definitions will be given later.)
In Section 3 we give an alternative proof 
a result in \cite{LM} on non-existence of hyperbolic leaves:

\begin{thm}[Lyubich-Minsky]\label{thm_feigenbaum_intro}
The regular part $\cR_f$ of a Feigenbaum map $f$ 
has only parabolic leaves.
\end{thm}

\parag{Cremer points and hedgehogs.}
Section 4 is devoted for rational functions 
with Cremer fixed points and their hedgehogs. 
For any small neighborhood of Cremer fixed point $\zeta_0$,
there exists an invariant continuum $H$ (a ``hedgehog")
containing $\zeta_0$, 
equipped with invertible ``sub-dynamics" $f|~H \to H$.

In Section 4, we will present a proof of 
the next result by following A.Ch\'eritat:

\begin{thm}[Lifted hedgehogs are irregular]\label{thm_Hhat_intro}
The invariant lift $\Hhat$ of $H$ is a continuum contained in the irregular part of the natural extension. 
\end{thm}

Since this natural extension has a continuum of irregular points,
one may expect to apply the classical Gross star theorem 
to find a hyperbolic leaf, as in \cite{KLR}.
However, the actual situation is not that good. 
It is still difficult to show the (non-)existence of
hyperbolic leaves without assuming the same conditions as \cite{KLR}.
Indeed, we will show that the irregular points in hedgehogs are 
not big enough to apply the Gross star theorem.

\parag{Acknowledgment.}
The authors would like to thank 
A.Ch\'eritat, J.Kahn, M.Lyubich, and L.Rempe for 
stimulating discussions on this subject.

\section{No hyperbolic leaves except Siegel disks}
Here we give a proof of Theorem \ref{thm_siegel}.
The idea is fairly simple and we will apply it to 
Feigenbaum quadratic maps in the next section.

Let us start with some terminologies:

\parag{Postcritical sets.}
Let $f:\C \to \C$ be a polynomial of degree more than one 
and we denote the set of its critical points by $C_f$.
The \textit{postcritical set} of $f$ is defined by 
$$
P_f :=  \overline{ \bigcup_{n \ge 1} f^n(C_f)}.  
$$
For example, if $f(z) = e^{2 \pi i \theta}z + z^2$ 
with $\theta$ of bounded type, then its Siegel disk $\Delta$ 
is a quasidisk and its boundary contains 
the only critical point $c_0 = -e^{2 \pi i \theta}/2$. 
Since the dynamics of $f:\overline{\Delta} \to \overline{\Delta}$ is 
topologically conjugate to an irrational rotation, 
we have $\partial \Delta = P_f$.

\parag{Deep points and uniform deepness.}
Let $K$ be a compact set in $\C$. 
For $x \in K$, let $\delta_x(r)$ denote
the radius of the largest open disk contained in $\D(x, r)-K$.
(When $\D(x,r) \subset K$, we define $\delta_x(r) := 0$.)
Then it is not difficult to check that the function
$(x,r) \mapsto \delta_x(r)$ is continuous.

We say $x \in K$ is a \textit{deep point} of $K$ 
if $\delta_x(r)/r \to 0$ as $r \to 0$.
\footnote{The original version due to C.McMullen 
\cite[\S 2.5]{Mc1} is as follows:
$x \in K$ is a \textit{deep point} if $\delta_x(r) = O(r^{1 + \al})$ for some $\alpha>0$. Here we mildly generalize it.}
For a subset $P$ of $K$, we say $P$ is \textit{uniformly deep} in $K$ if for any $\e>0$ there exists an $r_0$ such that 
for any $x \in P$ and $r<r_0$, we have $\delta_x(r)/r <\e$.

\parag{Deepness and measurable deepness.}
For a given measurable set $X \subset \C$, 
let $|X|$ denote its area.
For a compact set $K$ in $\C$ and an $X$ with $0<|X| < \infty$,
we define the {\it density} of $K$ in $X$ by
$$
\dens(K/X) := \frac{|K \cap X|}{|X|}.
$$
We say $x \in K$ is a {\it measurable deep point} of $K$ 
if $\dens(K/\D(x,r)) \to 1$ as $r \to 0$.
We say a subset $P$ of $K$ is {\it uniformly measurable deep} in $K$
if for any $\e>0$, there exists an $r_0 >0$ 
such that for any $x \in P$ and $r<r_0$ 
we have 
$$
1-\e < \dens(K/\D(x,r)) \le 1.
$$
\begin{lem}\label{lem_deepness}
Let $K$ be a compact set in $\C$. Then
\begin{enumerate}[\rm (1)]
\item
If $x \in K$ is a measurable deep point of $K$, 
then $x$ is a deep point.
\item 
If $P \subset K$ is uniformly measurable deep in $K$, 
then $P$ is uniformly deep.
\end{enumerate}
\end{lem}

\begin{pf}
It is enough to show (2):
Suppose that $P$ is uniformly measurable deep. 
Then for any $\e >0$, there exists an $r_0$ such that 
$1-\e <\dens(K/\D(x,r))$ for any $x \in P$ and $r<r_0$.
Since
$$
\dens(K/\D(x,r)) 
\le \frac{|\D(x,r)| -\pi \delta_x(r)^2}{|\D(x,r)|}
 =  1-\frac{\delta_x(r)^2}{r^2},
$$
we have $\delta_x(r)/r < \sqrt{\e}$. 
This implies that $P$ is uniformly deep. 
\QED\end{pf}

\parag{Uniform deepness of the postcritical set.}
Now let us go to the proof of Theorem \ref{thm_siegel}.
Set $f(z) = e^{2 \pi i \theta}z + z^2$, 
where $\theta$ is of bounded type.
Its Siegel disk is denoted by $\Delta$.

We will use the following result by C.McMullen \cite[\S 4]{Mc2}:

\begin{thm}[Uniform deepness of $P_f = \partial \Delta$]\label{thm_McMullen}
The postcritical set $P_f = \partial\Delta$ 
is uniformly measurable deep in $K_f$.
In particular, $P_f$ is uniformly deep by 
Lemma \ref{lem_deepness} above.
\end{thm}
Indeed, it is shown in \cite[\S 4]{Mc2} that 
there exist two positive constants $\alpha$ and $C$ such that for any $x \in \partial\Delta$, 
$
\dens(K_f/\D(x,r)) \ge 1-C r^{\alpha}
$
for sufficiently small $r>0$. 
\footnote{
The uniformity cannot be found in the statements of \cite{Mc2}, 
but it comes from the geometry of 
the Siegel disk $\Delta$ as a quasidisk. 
See also \cite[Cor.5, p.42]{BC} for a bit more elementary proof for the uniformity. }

\begin{pf}[Theorem \ref{thm_siegel}]
Let $\cR = \cR_f$ be the regular part of $\cN_f$, 
and $\Deltahat$ be the invariant lift of the Siegel disk $\Delta$.
We will show that for any leaf $L$ of $\cR-\Deltahat$
is parabolic. 

We first claim: 
\textit{Any leaf $L$ of $\cR-\Deltahat$ 
contains a backward orbit $\zhat$ that stays in
$\C- K_f$.}
In fact, if $L$ is parabolic, the projection 
$\pi:L \to \C$ has at most two exceptional values 
by Picard's theorem, and we can find such a $\zhat$.
If $L$ is hyperbolic and $\pi(L) \cap J_f \neq \emptyset$,
the claim is true because $\pi$ is an open map. 
In the remaining case 
$\pi(L)$ is contained in a Fatou component $U$. 
Since there is no critical point in the Fatou set, 
the projection $\pi:L \to U$ has to be a conformal isomorphism 
that cannot be extended.
This happens only when $L = \Deltahat$,
because $P_f = \partial \Delta$.
Hence the claim is justified.

\parag{Case 1.}
When $\zhat = \braces{z_{-n}}$ does not accumulate on $P_f =\partial\Delta$, the leaf $L = L(\zhat)$ is parabolic by a criterion of parabolicity by Lyubich and Minsky \cite[Cor.4.2]{LM}.

\parag{Case 2.}
Now let us assume that $\zhat = \braces{z_{-n}}$ accumulates on 
$P_f = \partial\Delta$.
Since $z_{-n}$ is contained in the basin at infinity, 
none of $z_{-n}$ hits the filled Julia set $K_f$,
in particular, none of $z_{-n}$ hits $\partial \Delta$ either.

By another criterion of parabolicity by Lyubich and Minsky \cite[Lem 4.4]{LM}, it is enough to show:
\begin{quote}
{\it
$\norm{Df^{-n}(z_0)} \to 0 ~(n \to \infty)$, 
where $Df^{-n}$ is the derivative of the branch of $f^{-n}$ sending $z_0$ to $z_{-n}$, 
and the norm is measured in the hyperbolic metric of $\C-\partial\Delta$.}
\end{quote}
Then $L = L(\zhat)$ is a parabolic leaf. 

Now set $\Omega:= \C - \overline{\Delta}$.
Then $z_{-n}$ is contained in $\Omega$ for all $n$.
Since $\Omega$ is topologically a punctured disk, 
it has a unique hyperbolic metric $\rho= \rho(z)|dz|$
induced by the metric $|dz|/(1-|z|^2)$ of constant curvature 
$-4$ on the unit disk.
To show the claim, it is enough to show 
$$
\norm{Df^{n}(z_{-n})}_\rho 
=
\frac{\rho(z_0)|Df^{n}(z_{-n})|}{\rho(z_{-n})}
\to \infty~~(n \to \infty),
$$
where the norm in the left is measured in the hyperbolic metric 
$\rho$. 

Here is a well-known property on $\rho$ that plays an important role
(See for example, \cite[Thm. 1-11]{Ah}):

\begin{lem}\label{lem_metric}
The hyperbolic metric $\rho_U = \rho_U(z)|dz|$ of any hyperbolic domain $U \subset \C$ 
is bounded by $1/d$-metric. 
More precisely,  
$$
\rho_U(z) \le \frac{1}{d(z, \partial U)}
$$
for any $z \in U$.
\end{lem}
We can check this by comparing $\rho_U$ and 
the hyperbolic metric on the disk 
of radius $r = d(z, \partial U)$ centered at $z$. 
Of course this lemma holds for $U = \Omega$.

Now it is enough to show:
\begin{equation}\label{eq:Df}
\norm{Df^{n}(z_{-n})}_\rho 
\asymp
\frac{|Df^{n}(z_{-n})|}{\rho(z_{-n})}
\ge
d(z_{-n}, \partial \Delta)|Df^{n}(z_{-n})|
\to \infty.
\end{equation}

Set $R_n: = d(z_{-n}, \partial \Delta)$.
By assumption, $R_n$ tends to $0$ by 
taking $n$ in a suitable subsequence.
Let $D_0$ denote the disk of radius $R_0$ centered at $z_{0}$, 
and let $U_n$ denote the connected component of $f^{-n}(D_0)$
containing $z_{-n}$. 
Since $D_0 \subset \Omega$, we have a univalent branch
$g_n:D_0 \to U_n$ of $f^{-n}$.
Set $v_n: = |Dg_n(z_0)| = |Df^n(z_{-n})|^{-1} >0$.
By the Koebe $1/4$ theorem, $g_n(D_0) = U_n$ contains the disk of radius $R_0 v_n/4$ centered at $z_{-n}$,  
and since $U_n \subset f^{-n}(\Omega) \subset \Omega$ we have
$R_0 v_n/4 \le R_n$.

\parag{Case 2-1:} 
First assume that $\liminf v_n/R_n = 0$.
If $n$ ranges over a suitable subsequence,
we have $v_n/R_n \to 0$ and thus (\ref{eq:Df}) holds. 
(Note that by the proof of \cite[Lem 4.4]{LM} we need (\ref{eq:Df}) 
only for a subsequence.)

\parag{Case 2-2:} 
Next consider the case when $\liminf v_n/R_n = q >0$.
We may assume that $n$ ranges over a subsequence 
with $\lim v_n/R_n =q$.

For $t>0$, let $tD_0$ denote the disk $\D(z_0, tR_0)$.
Since $D_0 = \D(z_0,R_0)$ is centered at a point in $\C-K$,
we can choose an $s<1$ such that $sD_0 \subset \C-K$.
By the Koebe 1/4-theorem,
$|g_n(sD_0)|$ contains $\D(z_{-n}, sR_0v_n/4) \subset \C-K$.

Let us take a point $x_n$ in $\partial\Delta$ 
such that $|x_n-z_{-n}| =R_n$. 
Then we have 
$$
\D(z_{-n}, sR_0v_n/4) \subset \D(x_n, 2R_n)
$$
and thus $\delta_{x_n}(2R_n) \ge sR_0v_n/4$.
Recall the assumption $v_n/R_n \sim q>0$ for $n \gg 0$.
This implies that the ratio 
$\delta_{x_n}(2R_n)/2R_n$ is bounded by a positive constant from below. 
However, $R_n =d(z_{-n}, \partial \Delta) \to 0$ by assumption and it contradicts to the uniform deepness of $P_f$ (Theorem \ref{thm_McMullen}). 
\QED
\end{pf}

\section{Feigenbaum quadratic polynomials}
In this section we give an alternative proof of 
Theorem \ref{thm_feigenbaum_intro} (\cite[Lem 4.6]{LM}).

\parag{Infinitely renormalizable quadratic maps.}
Let $U$ and $V$ be a topological disk with $U \Subset V \Subset \C$.
A proper branched covering $g:U \to V$ 
of degree two is called a \textit{quadratic-like map}.
A quadratic map $f_c$ is \textit{infinitely renormalizable}
if there is a sequence of quadratic-like maps 
$g_n:U_n \to V_n~(n = 0,1, \ldots)$ such that:
$0 \in U_n$ and it is the critical point of $g_n$;
$g_n$ is a restriction of $f_c^{p_n}$ on $U_n$
for some $p_n \in \N$;
and for each $n$ the ratio $p_{n + 1}/p_n$ is an integer $ \ge 2$.

\parag{Feigenbaum maps.}
An infinitely renormalizable quadratic map $f_c$
is called \textit{Feigenbaum} (or \textit{Feigenbaum-like}) if
there exist positive constants $b$ and $m$
independent of $n$ such that
$p_{n + 1}/p_n \le b$ and 
we can choose $U_n$ and $V_n$ satisfying
$\mathrm{mod}(V_n-\overline{U_n}) \ge m$. 

Now let us restate Theorem \ref{thm_feigenbaum_intro}:

\begin{thm}[Lyubich-Minsky]\label{thm_feigenbaum}
The regular part $\cR_f$ of a Feigenbaum quadratic map $f$ 
has only parabolic leaves.
\end{thm}

To apply the same argument as Theorem \ref{thm_siegel},
we first show:

\begin{prop}[Uniformly deep postcritical set]
\label{prop_feigenbaum}
The postcritical set $P_f$ of a Feigenbaum quadratic map $f$ is
uniformly deep in the filled Julia set $K_f$.
\end{prop}

\begin{pf}
By McMullen \cite[Thm.8.3]{Mc1}, 
each point in the postcritical set 
is a deep point of $K = K_f$, so we only need to show 
the uniform deepness.
Here we follow the idea of \cite[Cor. 5]{BC}.

Let us prove it by contradiction: 
Suppose that there exist an $\eta >0$, 
and sequences $x_i \in P_f$ and $r_i \searrow 0$ such that 
$\delta_{x_i}(r_i)/r_i \ge \eta$ for all $i \in \N$.
Then for each disk $D_i:= \D(x_i, r_i)$ 
we can find a disk $\Delta_i:= \D(y_i, \eta r_i)$
in $D_i -K$.
We may also assume that $x_i \to x \in P_f$ as $i \to \infty$.

For $n \in \N$ set 
$$
Q_n := \braces{z \in P_f \st 
\text{for any $r < 1/n$ }~\delta_{z}(r)/r \le \eta/20}.
$$
Then $Q_n$ is a closed set with $Q_{n} \subset Q_{n + 1}$
and $P_f = \bigcup Q_n$. 

It is known that the postcritical set $P_f$ is a Cantor set
with minimal invertible dynamics $f|\,P_f \to P_f$
\cite[Chap.8]{Mc1}.
In particular, $P_f$ itself is a complete metric space 
by the Euclidean metric on $\C$. 
Hence by the Baire category theorem there exists 
an $N$ such that $Q_N$ contains an open ball of $P_f$.
More precisely, there exists a round disk $V$ such that 
$E:= V \cap P_f \subset Q_N$. 

Next we take a neighborhood $U$ of $x$ 
and a univalent map $g:U \to g(U)$
such that $g(x) \in E$, $g(U \cap P_f) = g(U) \cap E$,
and $g(U \cap K) = g(U) \cap K$.
In fact, when $x$ does not land on the critical point,
 we just take a small disk $U$ around $x$ 
and $k \in \N$ such that $f^k(x) \in E$.
(Here we used the minimality of the dynamics on $P_f$.)
Then we have a univalent map $g = f^k|_U$ as desired.
Otherwise, there exist a $y \in E$ and a $k \in \N$ 
such that $f^k(y) = x$ (again by minimality) with 
$Df^k(y)\neq 0$. 
Then we take a small disk $U$ where 
univalent branch $g$ of such $f^k$ is defined.
For both cases, we may take $U$ as a disk of certain radius $r_0$
centered at $x$.

Now let us work with the deepness:
Let $v$ denote $|Dg(x)| \neq 0$ and 
$tU~(t >0)$ denote the disk $\D(x,tr_0)$.
By the Koebe distortion theorem, we can find a $t \in(0,1)$ 
such that 
$$
 \frac{v}{2} \le |Dg(z)| \le 2 v
$$
for any $z \in tU \subset U$.
When $i$ is sufficiently large, $D_i$ is contained in $tU$ and
we have 
$$
g(D_i) \subset \D(g(x_i), 2 v r_i)
~~\text{with}~~ g(x_i) \in E.
$$
Moreover, by the Koebe 1/4 theorem, 
$$
g(\Delta_i) \supset \D\paren{g(y_i), \frac{|Dg(y_i)| \eta r_i}{4}}
\supset \D\paren{g(y_i), \frac{v \eta r_i}{8}}.
$$
Since $g(\Delta_i) \subset \C -K$, we have
$
\delta_{g(x_i)}( 2 v r_i) \ge v \eta r_i/8
$
and thus 
$$
\frac{\delta_{g(x_i)}( 2 v r_i)}{2 v r_i} 
\ge \frac{\eta}{16} > \frac{\eta}{20}.
$$
This is a contradiction, since $g(x_i) \in E \subset Q_N$ and $r_i \to 0$. 
\QED
\end{pf}

\begin{pf}[Theorem \ref{thm_feigenbaum}]
We just follow the proof of Theorem \ref{thm_siegel}.
Take a leaf $L$ of $\cR_f$.
To show $L$ is parabolic we only need to check: 
\begin{enumerate}[(i)]
\item
We can always take a backward orbit $\zhat \in L$ contained in $\C-K_f$.
\item 
The hyperbolic metric on $U = \C-P_f$ is hyperbolic in order 
to apply Lemma \ref{lem_metric}.
\item 
The postcritical set $P_f$ is uniformly deep in $K_f$.
\end{enumerate}
(i) is clear since $K_f$ has no interior. 
(ii) is also clear because $P_f$ contains at least two points.
(iii) is checked by Proposition \ref{prop_feigenbaum}.
\QED
\end{pf}

\parag{Remark.}
As Lyubich and Minsky's original proof of Theorem
\ref{thm_feigenbaum_intro} in \cite[Lem. 4.6]{LM}, 
it would be possible to show Theorem \ref{thm_siegel} 
without using the deepness of the postcritical set.

The virtue of our proofs presented here 
is that we only concern with the geometry of the Julia set 
and the postcritical set, 
and there is no need to look at the dynamics in detail.
One may expect to apply it to the Cremer case. 
But so far the authors do not have any example of 
Cremer quadratic map with uniformly deep postcritical set.

\section{Lifted hedgehogs and the Gross criterion}

Let $f:\Cbar \to \Cbar$ be 
a rational function of degree more than one 
that has a Cremer fixed point $\zeta_0$. 
Here is a fundamental result due to 
P\'erez-Marco (See \cite[Thm. 1]{PM1} and \cite[Thm. 2]{PM2}):

\begin{thm}[Cremer hedgehogs]\label{thm_PM}
Let $D_0$ be any neighborhood of $\zeta_0$ where $f$ is univalent and there exists a univalent branch $g:D_0 \to \C$ of $f^{-1}$ with $g(\zeta_0) = \zeta_0$. Then for any Jordan neighborhood $D \Subset D_0$ of $\zeta_0$, there exists a compact set $H \subset \overline{D}$ containing $\zeta_0$ with the following properties:
\begin{enumerate}
\item
$H$ is a connected and full continuum in the Julia set.
\item
$g(H) = H$, $f(H) = H$, and $H \cap \partial D \neq \emptyset$.
\item
There exists a subsequence $\braces{n_k} \subset \N$ such that $f^{n_k}|_H \to \id~~(k \to \infty)$.
\end{enumerate}
\end{thm}

We say $H$ a \textit{hedgehog} of the Cremer point $\zeta_0$.
By this theorem we have an invertible ``sub-dynamics" $f:H \to H$. 
Let $\cN = \cN_f$ be the natural extension of $f:\Cbar \to \Cbar$ and $\Hhat$ be the {\it lifted hedgehog} in $\cN$, that is, 
the set of backward orbits in remaining in $H$.

Here we restate Theorem \ref{thm_Hhat_intro}:
\begin{thm}[Lifted hedgehogs are irregular]\label{thm_Hhat}
Any point $\zhat$ in $\Hhat$ is irregular in $\cN$. In particular,
$\Hhat$ forms a continuum of irregular sets.
\end{thm}
The authors learned the principal idea of the proof 
by A.Ch\'eritat.

\begin{pf}
Suppose that $\zhat = (z_0,z_{-1}, \ldots) \in \Hhat$ is a regular point, that is, there exists a neighborhood $U_0$ of $z_0 \in H$ such that its pull-back along $\zhat$ is eventually univalent. 

By replacing $z_0$ by  $z_{-n}$ with sufficiently large $n$, we may assume that whole pull-back $\cdots \to U_{-1} \to U_0$ along $\zhat$ is univalent. 
Then $g^n|~ U_0 \to \Cbar$ makes sense for all $n$ and $U_{-n} = g^n(U_0)$ does not touch the critical set $C_f$ and its preimage $f^{-1}(C_f)$. Since $f$ has degree more than one, the family $\braces{g^n|U_0}$ is normal and thus it has a holomorphic sequential limit $\lim g^{n_k} = G:U_0 \to \Cbar$.

If $G$ is not constant, then the set $G(U_0)$ is an open set intersecting with $H$. Since $f^{n_k} \circ G \to \id~(k \to \infty)$, $f^{n_k}$ is normal on $G(U_0)$. This implies that $G(U_0)$ is contained in the Fatou set, but it contradicts the fact that $H$ is contained in the Julia set.

Now $G$ must be constant. However, we may also assume that $g^{n_k}|{U_0 \cap H} \to \id$ by the theorem above. This is again a contradiction.

Since the backward action $g:H \to H$ is a homeomorphism, $H$ and $\Hhat$ is homeomorphic. Hence $\Hhat$ is a compact continuum. 
\QED
\end{pf}

\subsection{Does a hedgehog generate any hyperbolic leaves?}
By the theorem above, we have a fairly big set of irregular points
in the natural extension generated by the Cremer hedgehogs. 
Hence, according to \cite[Lem.3.3]{KLR}, 
one would expect to apply the Gross star theorem
to find hyperbolic leaves. 
However, we will see that these hedgehogs are not big enough 
to apply the Gross star theorem.

For simplicity we consider 
the case when $f(z) = z^2 + c$ is a quadratic polynomial. 
Let $P$ and $J$ denote the postcritical set and the Julia set.
(Conventionally we remove $\infty$ from quadratic postcritical sets.)
For the natural extension $\cN = \cN_f$, 
let $\pi_{-n}:\cN \to \Cbar$ denote the projection 
$\pi_{-n}(\zhat) = z_{-n}$, the $n$-th entry of $\zhat$.
Instead of $\pi_0$ we use the notation $\pi$ as before.

\parag{The Gross criterion.}
Fix any $z_0 \in \C-P$. 
Then each $\hat{z} \in \pi^{-1}(z_0)$ is regular in $\cN$. 
In particular, the projection $\pi:L(\zhat) \to \Cbar$ is locally univalent near $\pi:\zhat \mapsto z_0$.

Let $\ell(\theta)~(\theta \in [0,2\pi))$ denote the half-line given by $\ell(\theta) := \braces{z_0 + re^{i\theta} \st r \ge 0}$. 
By using the Gross star theorem, \cite[Lem.3.3]{KLR} claims:
\textit{
if $L(\zhat)$ is isomorphic to $\C$, 
then for almost every angle $\theta \in [0,2\pi)$ 
the locally univalent inverse $\pi^{-1}: z_0 \mapsto \zhat$ 
has an analytic continuation along the whole half-line $\ell(\theta)$.}
Hence to show $L(\zhat)$ is hyperbolic, we should show:
\begin{quote}
$(\ast)$: \textit{There exist a $\zhat \in \pi^{-1}(z_0)$ 
and a set of $\theta \in [0, 2\pi)$ of positive length
such that the analytic continuation of $\pi^{-1}: z_0 \mapsto \zhat$ 
along $\ell(\theta)$ hits an irregular point at some $z = z_0 + re^{i \theta}$.}
\end{quote}
This is what they do in \cite[Prop.3.2]{KLR} to find hyperbolic leaves that intersect the Julia set.

\parag{Where the irregular points come from?}
Now we assume that $(\ast)$ holds in our setting.
We will claim that for almost every $\theta$, 
the irregular point corresponding to 
the half-line $\ell(\theta)$ 
does not belong to the lifted hedgehog $\Hhat$:
that is, the hyperbolic leaf is not generated 
by $\Hhat$.

Since the critical orbit $CO:= \braces{f^n(0)}_{n \ge 1}$ 
is a countable set, 
the set of angles $\theta$ with $\ell(\theta) \cap CO \neq \emptyset$
is at most countable, i.e., a null set of angles. 
We may forget about such angles for our purpose.
So we define the set $\Theta \subset [0,2\pi)$ as follows: 
$\theta \in \Theta$ if $\ell(\theta) \cap CO = \emptyset$
and $\ell(\theta) \cap H \neq \emptyset$.
Since $H$ is a continuum, $\Theta$ is a set of angles with positive length.

One more fact we need to recall is that the hedgehog $H$ does not intersect with $CO$. (But $H$ is contained in $P = \overline{CO}$. See \cite{Ch}.) In fact, if $H \cap CO \neq \emptyset$, the only critical point of $f(z) = z^2 + c$ eventually captured in $H ( = f^n(H)$ for all $n \ge 0)$. 
However, this contradicts the fact that larger hedgehogs are also contained in $P$. 

Hence hedgehog $H$ has a ``backward orbit" $\cdots \to H_{-2} \to H_{-1} \to H$ with $H_{-n}$ homeomorphic to $H$. In particular, $\pi^{-1}(H)$ consists of homeomorphic copies of $H$.

\begin{prop}\label{prop_}
For any $\ell(\theta)$ with $\theta \in \Theta$, there exists a unique path $\ell(\theta,H)$ in the natural extension such that 
it intersects with $\Hhat$ and 
$\pi: \ell(\theta,H) \to \ell(\theta)$ is a homeomorphism. 
\end{prop}

\begin{pf}
Since $\ell(\theta)$ does not contain the points of $CO$, 
its preimage by $f^n$ consists of exactly $2^n$ disjoint curves. 
We can uniquely choose one of them that touches $H$,
since $\ell(\theta) \cap H \neq \emptyset$, $H \cap CO =\emptyset$, and $f\mid H \to H$ is a homeomorphism. 
Now the backward orbits remaining in such a curve for each $n$ form a unique path in the natural extension that has a unique initial point $\zhat(\theta,H) \in \pi^{-1}(z_0)$ and passes though a point in $\Hhat$. This is the desired path $\ell(\theta,H)$. 
\QED
\end{pf}

\parag{Remark.}
Note that the analytic continuation along $\ell(\theta)$ of the germ $\pi^{-1}$ near $\pi:\zhat(\theta,H) \mapsto z_0$ will have singularity at least when it hits $\Hhat$.

Next we prove:
\begin{prop}\label{prop_injectivity}
The map $\chi:\Theta \to \pi^{-1}(z_0)$ defined by 
$\chi(\theta):= \zhat(\theta,H)$ is injective.
\end{prop}
Recall that every $\theta \in \Theta$ uniquely determines 
$\ell(\theta,H)$ and $\zhat(\theta, H)$. 
By the proposition above,
if there is a $\zhat = \zhat(\theta, H)$ as in $(\ast)$, 
then for any $\theta' \in \Theta-\{\theta\}$,  
the analytic continuation of 
the locally univalent inverse 
$\pi^{-1}: z_0 \mapsto \zhat(\theta, H)$ 
along $\ell(\theta')$ does not hit $\Hhat$.
Hence it can only hit an irregular point 
that is \textit{not} contained in $\Hhat$.

\begin{pf}
Let us take $\theta \in \Theta$ and $\theta' \in \Theta$ 
with $\theta \neq \theta'$. 
Then within the bounded region enclosed by $\ell(\theta)\cup \ell(\theta') \cup H$, we have at least one point in $CO$. (Actually infinitely many.)  
Let $z_{-n}(\theta,H) \in f^{-n}(z_0)$ denote $\pi_{-n}(\zhat(\theta,H))$, i.e., the starting point of the $n$-th pull-back of $\ell(\theta)$ along $\Hhat$. 
Now it is enough to claim:
\textit{
There exists an $n$ such that $z_{-n}(\theta',H)$ 
is different from $z_{-n}(\theta,H)$. }

Let $U$ be the connected component of the region enclosed by $\ell(\theta)\cup \ell(\theta') \cup H$ containing $z_0$ on its boundary.
Suppose that $z_{-n}(\theta,H) = z_{-n}(\theta',H) =: z_{-n}$ for all $n$. 
Then we have a pull-back $\braces{U_{-n}}$ of $U = U_0$ 
that has $\braces{z_{-n}}$ on its boundary. 

If $U_{-n}$ contains a critical point, 
$z_{-n}(\theta',H)$ and $z_{-n}(\theta,H)$ can not be the same.
 (See Figure \ref{fig_H}.) 
\begin{figure}[htbp]
\centering
\includegraphics[width=.55\textwidth]{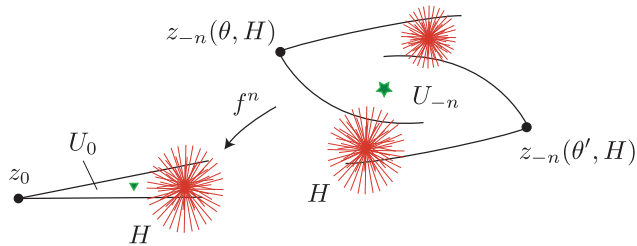}
\caption{
The star is the critical point
and it is mapped to the triangle by $f^n$.}\label{fig_H}
\end{figure}
Thus $U_{-n}$ does not contain the critical point for all $n$,
and we have a univalent branch $h_n:U \to U_{-n}$ of $f^{-n}$.
Since $h_n(U)$ avoids $0$, $c$, and $\infty$, 
$\braces{h_n}$ forms a normal family. 
By a similar argument as in the theorem above, 
any sequential limit cannot be an open map and 
they must be constant.
However, this is impossible, because the compact set $\partial U \cap H$ contains a continuum of definite diameter, and the action 
on $f^{n_k}|\partial U \cap H$ will be arbitrarily close to identity. \QED
\end{pf}

In conclusion, the lifted hedgehog $\Hhat$ is not big enough to apply the Gross criterion. On the other hand, according to the technique of 
Theorem \ref{thm_siegel}, it seems reasonable to  

\parag{Conjecture.}
\textit{There exists a Cremer quadratic 
polynomial whose regular part has no hyperbolic leaf.}

\medskip

{\small
\noindent 
Carlos Cabrera. 
Instituto de Matem\'{a}ticas,
Unidad Cuernavaca, UNAM. 
Cuernavaca, Mexico.

\medskip

\noindent 
Tomoki Kawahira. 
Graduate School of Mathematics, 
Nagoya University. 
Nagoya, Japan. 
}

\begin{thebibliography}{10}
\bibitem[Ah]{Ah} L.V. Ahrlfors. 
\textit{Conformal Invariants}. 
McGraw-Hill, 1973.

\bibitem[BC]{BC}
X. Buff and A. Ch{\'e}ritat. 
Quadratic Julia sets with positive area.
\textit{Preprint}. ({\tt arXiv:math/0605514v2}).

\bibitem[Ch]{Ch}
D.K. Childers. Are there critical points on the boundaries of mother hedgehogs?  
\textit{Holomorphic dynamics and renormalization},   
Fields Inst. Commun., {\bf 53} (2008), 75--87.

\bibitem[KL]{KL} V.A. Kaimanovich and M. Lyubich. 
Conformal and harmonic measures on laminations associated with rational maps.
\textit{Mem. Am. Math. Soc.}, {\bf 820}, 2005. 

\bibitem[KLR]{KLR} J. Kahn, M. Lyubich, and L. Rempe.  
A note on hyperbolic leaves and wild laminations of rational functions.
\textit{J. Difference Equ. Appl.}, {\bf 16} (2010), no. 5--6, 655--665. 

\bibitem[LM]{LM} M. Lyubich and Y. Minsky.
Laminations in holomorphic dynamics.
\textit{J. Diff. Geom.} {\bf 49} (1997), 17--94.

\bibitem[Mc1]{Mc1} C. McMullen. 
\textit{Renormalization and 3-manifolds which fiber over the circle}. Annals of Math Studies 142, Princeton University Press, 1996.

\bibitem[Mc2]{Mc2} C. McMullen. 
Self-similarity of Siegel disks and Hausdorff dimension of Julia sets. {\it Acta Math.} {\bf 180} (1998), no.2, 247--292.

\bibitem[PM1]{PM1} R. P\'erez-Marco. 
Fixed points and circle maps.
{\it Acta Math.} 
{\bf 179} (1997), 243--294.

\bibitem[PM2]{PM2} R. P\'erez-Marco. 
Sur une question de Dulac et Fatou.
{\it C. R. Acad. Sci. Paris S\'er. I Math.} 
{\bf 321} (1995), no. 8, 1045--1048.

\end{thebibliography}
\end{document}